\input amstex
\magnification=1200
\documentstyle{amsppt}
\NoRunningHeads
\TagsOnRight
\baselineskip 20pt
\pageheight{8.7truein}
\pagewidth{5.5truein}

\NoBlackBoxes
\define\La{\Lambda}
\define\la{\lambda}
\define\BLa{\boldsymbol\Lambda}
\define\ga{\gamma}
\define\CA{{\Cal A}}
\define\CB{{\Cal B}}
\define\CF{{\Cal F}}
\define\CM{{\Cal M}}
\define\CC{{\Bbb C}}
\define\RR{{\Bbb R}}
\define\ZZ{{\Bbb Z}}
\define\({\lf(}
\define\){\rt)}
\define\[{\lf[}
\define\]{\rt]}
\define\lf{\left}
\define\rt{\right}


\topmatter
\title A perturbative SU(3) Casson invariant
\endtitle
\author S.E. Cappell, R. Lee, E.Y. Miller\endauthor
\footnote""{During the research of this paper, all
three authors were supported by grants from the National Science Foundation.}

\endtopmatter

\bigskip

\subhead 1.  Introduction \endsubhead

    From a gauge theory viewpoint, the well-known
$SU(2)$-Casson invariant $\lambda_{SU(2)}(X)$ of an integral homology
3-sphere $X$ can be regarded as 
the number, counted with sign, of flat $SU(2)$-connections on $X$ after
making a suitable perturbation of the curvature equation \cite{T}. 
In Casson's original treatment, $\lambda_{SU(2)}$ was
obtained from a finite dimensional, symplectic setting, as
the intersection number in a representation varity of two 
perturbed Lagrangian subvarieties 
associated to a Heegaard decomposition of $X$ (see \cite{AM}).
In both these
gauge-theoretic and symplectic settings, 
the fact that pertubations were used
in the definition and that large scale 
perturbations are permissible underlay
remarkable properties of the Casson invariant, such as surgery formulae.
In this paper we solve
the problem of defining a (fully) perturbative $SU(3)$ generalization, 
$\Lambda_{SU(3)}(X)$,
of the Casson invariant, and begin the study of its properties. Some of these
recall well-known facts about the $SU(2)$-Casson invariant:

(1) An integrality property: $4\cdot\Lambda_{SU(3)}(X) \in \Bbb Z$.

(2) In the cases computed here, for $1/k$-surgery on some torus knots, the
invariants are given by  quadratic polynomials in $k$, for $k$ positive 
(resp. negative) while
in the $SU(2)$ case they are linear.

(3) It is preserved under the change of orientation, just as the $SU(2)$-
invariant is reversed.

On the other hand, it differs intriguingly from the $SU(2)$-invariant in
that the polynomials giving the values for $1/k$-surgery on the torus knots
in (2) for $k$ positive are not the same as those for $k$ negative.

Our investigation has benefited greatly from the excellent series of 
recent articles of Boden-Herald and of Boden-Herald-Kirk-Klassen \cite{BH 1,2},
\cite{BHKK}. In \cite{BH 1} 
a different gauge-theoretic generalization,
$\lambda_{SU(3)}(X)$, of the Casson invariant to $SU(3)$ was introduced 
using - and allowing - only small perturbations; it is thus not fully 
perturbative. Among the important properties Boden-Herald obtained for 
their invariant are: $\lambda_{SU(3)}$ is independent of orientation,
$\lambda_{SU(3)}(X)=\lambda_{SU(3)}(-X)$, and has a connect sum formula
$\lambda_{SU(3)}(X_1\#X_2)=\lambda_{SU(3)}(X_1)+\lambda_{SU(3)}(X_2)+
4\lambda_{SU(2)}(X_1)\lambda_{SU(2)}(X_2)$ (see \cite{BH 1,2}). 
In the paper \cite{BHKK}, there
are impressive calculations of this invariant for $1/k$-surgery
on some torus knots, with the result that the values are given by various 
rational functions in $k$, 
cubic polynomials divided by linear polynomials in their
cases. As is already evident from their calculations,
in special cases 
$\lambda_{SU(3)}$ takes values which are fractions 
with varying denominators; moreover this  would follow more generally from 
a conjecture on Chern-Simons invariants.
Thus this contrasts with the integrability 
property of the invariant considered here.

Some years ago, in \cite{CLM}, we proposed a program for defining a 
generalized $SU(n)$-Casson invariant based on a Lagrangian intersection
number of perturbed subvarieties in 
the $SU(n)$-representations of $\pi_1(X)$.  
That program proposed using in such a definition, correction terms 
obtained from combinations of 
tangential and normal Maslov indices along the singular strata of
reducible representations. In part to understand these correction terms, 
we studied the relation between Maslov index and spectral flow in 
\cite{CLM 1,2} and the different definitions of $SU(2)$-Casson invariants
for rational homology spheres in \cite{CLM 3}. The present effort 
could be viewed as a modification and completion of the program of 
\cite{CLM} for $SU(3)$. The new ingredient in the definition is a 
further term which involves the boundary maps of the mod-$2$ Floer
chain complex \cite{F}. 
It is this extra term which makes the invariant well-defined
and, as in Theorem (3.4), fully perturbative as we had wished.

We now provide a precise comparison of these two invariants 
$\lambda_{SU(3)}(X)$
and $\Lambda_{SU(3)}(X)$. Recall that the $SU(2)$-Casson invariant, 
$\lambda_{SU(2)}(X)$, was reformulated by Taubes \cite{T} in a 
gauge-theoretic setting, as the sum:
$$
\lambda_{SU(2)}(X) = (-1)\sum \Sb {[A]\in \Cal M^*_{SU(2),h}}\endSb 
(-1)^{SF(\theta,A,h;su(2))} \tag1.1
$$
where $[A]$ runs through all gauge equivalent classes of $h$-perturbed
$SU(2)$-connections. 
\footnote{The first minus sign (-1) is explained in \cite{KK} also page 5 of
\cite{BH}.}
The sign $(-1)^{SF(\theta,A,h;su(2))}$ is specified by the spectral flow 
$SF(\theta,A,h;su(2))$ associated to a path of connections from the trivial
connection $\theta$ to $A$. In theory, this last spectral flow depends on the
choice of paths; however the ambiquity equals to $0$ (mod $8$) and 
thus disappears 
when we form the sign $(-1)^{SF(\theta,A,h;su(2))}$.

In the work of Boden-Herald \cite{BH}, 
as briefly reviewed in \S2 below,
the invariant $\lambda_{SU(3)}(X)$ is given by
$$
\align
\lambda _{SU(3)}(X) &=\lambda^{'}_{SU(3)}(X)+\lambda^{''}_{SU(3)}(X) \tag1.2 \\
\lambda^{'}_{SU(3)}(X) &=\sum \Sb {[A]\in \Cal M^*_{SU(3),h}}\endSb
(-1)^{SF(\Theta,A,h;su(3))} \\
\lambda^{''}_{SU(3)}(X) &=\sum \Sb {[A]\in \Cal M^*_{SU(2),h}}\endSb  
(-1)^{SF(\theta,A,h;su(2))}[SF(\theta,A,h;\Bbb C^2)-2cs(\Hat A)+1]
\endalign
$$
after making a ``small'' perturbation $h$. The correction term 
$\lambda^{''}_{SU(3)}$ is introduced because the number 
$\lambda^{'}_{SU(3)}$ of $h$-perturbed flat, irreducible, $SU(3)$-
connections in $\Cal M^*_{SU(3),h}$ depends on the choice of 
perturbations. Given two perturbations $h_0, h_1$, we can 
connect them up by a family of small perturbations $h_t,0\leq t\leq 1$.
Along this path, there would exist a cobordism joining points in
$\Cal M^*_{SU(3),h_0}$ and $\Cal M^*_{SU(3),h_1}$ but for the 
phenomena of irreducible $SU(3)$-connections sinking into or
emerging from the $SU(2)$-stratum. Whenever this occurs, a corresponding  
integer jump occurs in the normal spectral flow 
$SF(\theta,A,h;\Bbb C^2)$. Thus the discrepancy in $\lambda^{'}_{SU(3)}(X)$
is compansated by the sum 
$\Sigma (-1)^{SF(\theta,A,h;su(2))}[SF(\theta,A,h;\Bbb C^2)]$.

However, the above spectral flow $SF(\theta,A,h;\Bbb C^2)$
depends on the choice of paths from the trivial
connection $\theta$ to $A$. By definition, a ``small'' perturbation
has the property that the $h$-perturbed flat, irreducible, 
$SU(2)$-connections $[A]\in \Cal M^*_{SU(2),h}$ is within 
$\epsilon$-distance of a unique component $\hat A$ in 
the space $\Cal M^*_{SU(2)}$ of flat connections. In particular, we have
a well-defined path class 
$\alpha$ from $A$ to an element in the component $\Hat A$.
Given such a component $\Hat A$, we can also choose a path $\beta$ 
connecting
an element in $\Hat A$ to the trivial connection $\theta$ and using this 
path we can calculate the Chern-Simons invariant $cs(\hat A)$. On 
the other hand, the composite $\beta \circ \alpha$ provides a way to
connect up $A$ with $\theta$, and hence a spectral flow invariant 
$SF(\theta,A,h;\Bbb C^2)$. Although both $cs(\Hat A)$
and $SF(\theta,A,h;\Bbb C^2)$ depend on the choice of the path $\beta$, 
the ambiguities cancel each other and the combination 
yields a well-defined term
$[SF(\theta,A,h;\Bbb C^2)-2cs(\Hat A)+1]$ in $(1.2)$.

Now the present perturbative $SU(3)$-Casson invariant $\Lambda_{SU(3)}(X)$
is given by the formula:
$$
\align
\Lambda _{SU(3)}(X) &=\Lambda^{'}_{SU(3)}+\Lambda^{''}_{SU(3)}(X) 
-(1/4){\text{Floer}}(X,h) \tag1.3 \\
\Lambda^{'}_{SU(3)}(X) &=\sum \Sb {[A]\in \Cal M^*_{SU(3),h}}\endSb
(-1)^{SF(\Theta,A,h;su(3))} \\
\Lambda^{''}_{SU(3)}(X) 
&=\sum \Sb {[A]\in \Cal M^*_{S(U(1)\times U(2)),h}}\endSb
(-1)^{SF(\theta,A,h;s(u(1)\times u(2)))}[SF(\theta,A,h;\Bbb C^2) \\
&-(1/4)\,SF(\theta,A,h;s(u(1)\times u(2)))+5/8] \\
{\text{Floer}}(X,h) &=\sum_{p=0}^7 (-1)^pdim_{Z/2}({\text{Image d:}}
\quad FC_{p+1}(X,h)\to FC_{p}(X,h))
\endalign
$$
Here the first term $\Lambda^{'}_{SU(3)}(X)$ is the same as 
$\lambda^{'}_{SU(3)}(X)$. In the second term $\Lambda^{''}_{SU(3)}(X)$,
the normal spectral flow $SF(\theta,A,h;\Bbb C^2)$ is the same as that in
$\lambda^{''}_{SU(3)}(X)$ while the Chern-Simons term $cs(\Hat A)$ is replaced 
by 1/4 of the tangential spectral flow,
$(1/4)SF(\theta,A,h;s(u(1)\times u(2)))$.
The combination $[SF(\theta,A,h;\Bbb C^2)-
(1/4)SF(\theta,A,h;s(u(1)\times u(2)))]$ was shown in \cite{CLM} to be 
independent of the choice of paths connecting $\theta$ to $[A]$ and
has the advantage of being free from the 
restrictive assumption of small perturbations. 
Unfortunately the tangential spectral flow $SF(\theta,A,h;s(u(1)\times u(2)))$
also creates a problem of its own. For a family of perturbations $h_t$,
a pair $(A_t(1), A_t(2))$
of $h_t$-perturbed flat, irreducible, $SU(2)$-connections 
can be created or destroyed through their collision at
a birth-death point (the analogue of Whitney disk cancellation in 
the context of 
finite dimensional handle decompositions). Whenever this happens, 
the terms in the sum
$\Sigma SF(\theta,A,h;s(u(1)\times u(2)))$ corresponding to
$(A_t(1), A_t(2))$ 
will cause a jump and so
$\Lambda^{'}_{SU(3)}+\Lambda^{''}_{SU(3)}(X)$ is not a well-defined invariant.

Analogues of such problems of jumps have been 
studied in parametrized Morse theory, but here we have to adjust this to 
the infinite dimensional gauge space
with the Chern-Simons functional as the Morse function.
Although the Floer homology $FH_{*}(X)$ with $Z/2$-coefficients
\footnote{We can also work with Floer homology in integer or 
other coefficients.}
is well-defined, its Floer
chain groups $FC_{*}(X,h)$ varies precisely because of
the existence of these birth-death points. Indeed,
a fixed integer jump occurs in ${\text{Floer}}(X,h_t)$ when
$h_t$ goes through such a 
birth or death point. Hence ${\text{Floer}}(X,h)$ can be used as
a correction term for the discrepancy in 
$\Sigma SF(\theta,A,h;s(u(1)\times u(2)))$. 
Detailed analysis of $\Lambda^{'}_{SU(3)(X)},\Lambda^{''}_{SU(3)(X)},
\text{Floer}(X,h)$ as well as the proof that $\Lambda_{SU(3)}(X)$ is 
well-defined (Theorem 3.4) can be found in \S3.

Despite the differences between $\lambda_{SU(3)}$ and $\Lambda_{SU(3)}$,
they also share some properties. 
For example, they are independent of orientation
(see Proposition 4.5 for $\Lambda_{SU(3)}$) and have connect sum formulae.
Due to the Floer correction term, the formula for $\Lambda_{SU(3)}$
is more complicated than its counterpart in \cite{BH 2}, as it involves
the Floer chain complex of the connected sum which is a subtle aspect
of Floer homology theory (see \cite{Fu},\cite{Li}). 
The proof of this connect
sum formula for $\Lambda_{SU(3)}$ is in \S4.

In \S5, we provide explicit calculations of our invariant for the
Brieskorn spheres $\Sigma(2,q,2qk\pm 1), q=3,5,7,9,$ which can also
be obtained from $\mp 1/k$ surgery on $(2,q)$-torus knots. Our 
results are parallel to those in \cite{BHKK} where 
$\lambda_{SU(3)}(\Sigma(2,q,2qk\pm 1))$ in the same range are 
computed. However we have to calculate the 
spectral flow $SF(\theta,A,h; s(u(1)\times u(2)))$ for all flat, irreducible
$SU(2)$-connections $[A]$. In \cite{FS}, Fintushel-Stern calculated these
spectral flows and their results are tailor-made for us (see Theorem 5.1).

As mentioned before, Casson's $SU(2)$-invariant was first defined
using Heegaard decomposition and intersection of perturbed Lagrangians in
the representation varieties. We briefly discuss how the 
representation-theoretic analogue of the present gauge-theoretic treatment of
$\Lambda_{SU(3)}$ would proceed, 
as this was the context envisioned in \cite{CLM}: 
Using a Heegaard decomposition, we can write $X$ as a union
$X_1\cup X_2$ of two handle bodies $X_1, X_2$ glued along a Riemann surface
$\Sigma$. Then the moduli space $\Cal M_{SU(3)}(X)$ of flat $SU(3)$-connections
can be identified with the intersection of the Lagrangian subspaces
$R_{SU(2)}(X_i)=Hom(\pi_1(X_i),SU(3))/SU(3)$ inside 
$R_{SU(3)}(\Sigma)=Hom(\pi_1(\Sigma),SU(3))/SU(3)$. 
After a suitable Hamiltonian perturbation, the Maslov indices
at the reducibles are defined and a Floer correction introduced.
Then the symplectic definition of $\Lambda_{SU(3)}$ is the same
as in (1.3). Indeed
to define $\text{Floer}(X,h)$,
it is natural to consider a symplectic Floer homology theory based on the
intersection of $R_{SU(2)}(X_i)$ in the $SU(2)$-stratum $R_{SU(2)}(\Sigma)$.
In this direction, there are the work of Lee-Li \cite{LL}
which treats the singular nature of $R_{SU(2)}(\Sigma)$ and the  
work of Sullivan \cite{S} 
which addresses the change of Floer chain complexes in
the smooth context under perturbations.

Finally, the general methodology introduced here to define fully perturbative
invariants using $\text{Floer}(X,h)$ may appear complicated in that 
this term has a ``tertiary'' character, being the correction to the Maslov 
index correction term along singularities. But this method opens up for 
$\Lambda_{SU(3)}$, and perhaps much more generally, 
the possibility of intriguing relations with still unknown Floer theories. 
In particular, as $8\cdot\Lambda_{SU(3)}$ 
\footnote{Although $4\cdot\Lambda_{SU(3)}$ is an integer, it is more natural 
to consider $8\cdot\Lambda_{SU(3)}$ as an Euler characteristic.}
is an integer invariant, it suggests 
the existence of a $SU(3)$-Floer homology with $8\cdot\Lambda_{SU(3)}$
as its Euler characteristic.  

\vskip .25in

\subhead \S2. Review of the work of Boden, Herald, Kirk and Klassen\endsubhead

Let $X$ be an oriented, integral homology $3$-sphere and let $\Cal A$ be the 
space of smooth, 
$SU(3)$-connections on the trivial product bundle $P = X\times 
SU(3)$.  This last space $\Cal A$ is an infinite dimensional affine space and 
in fact by fixing a trivial product connection 
$\theta$ on $P$, we can identify 
$\Cal A$ with the space $\Omega^1(X, Ad P) = \Omega^1(X, su(3))$ of 
$su(3)$-valued $1$-form on $X$.

Let $\Cal G = \text{Map}(X, SU(3)) = C^\infty (X,SU(3))$ denote the gauge group of 
$SU(3)$-bundle automorphisms $g:P\to P$ of $P$.  Then as these gauge 
transformations change the bundle structure and hence the connections 
$A\to g\cdot A = gAg^{-1} + gdg^{-1}$, 
they give rise to an action of $\Cal G$ on $\Cal A$ with 
$\Cal B =\Cal A/\Cal G$ as quotient.  
This action is not free, and according to the 
isotropy subgroup there is the natural Whitney stratification on $\Cal A$ and 
also on the orbit space $\Cal B = \Cal A/\Cal G$.  A $SU(3)$-connection $A$ in 
$\Cal A$ is said to be irreducible if its isotropy 
subgroup consists of constant maps 
to $\bold Z(SU(3)) = \bold Z/3$.  Altogether these irreducibles form the top 
stratum $\Cal A^*$ and its quotient 
$\Cal B^* = \Cal A^*/\Cal G$ has a structure of 
pre-Banach manifold.

Below the top stratum, there are strata whose isotropy subgroups are 
respectively $U(1), S(U(1)\times U(1)\times U(1))$, 
$S(U(1)\times U(2))$ and $SU(3)$,  They corresond 
to the situation where the underlying $3$-dimensional complex vector bundles 
and connections are decomposed into:

\flushpar{\bf (2.1)}
\roster
\item"(a)" A sum $L\oplus Q$ of line bundle $L$ and a $2$-plane bundle 
$Q$ with structure group $S(U(1)\times U(2))$.
\item"(b)" A sum $L_1\oplus L_2\oplus L_3$ of three line bundles 
$L_1,L_2,L_3$ which are all different 
$L_1\ne L_2 \ne L_3$ and with structure group 
$S(U(1)\times U(1)\times U(1))$.
\item"(c)"  A sum $L_1\oplus L_2\oplus L_3$ of three line bundles, 
two of which are the same and with structure group $S(U(1)\times U(1))$.
\item"(d)"  A sum $L_1\oplus L_1\oplus L_1$ 
of three isomorphic line bundles with structure group $\bold Z/3$.
\endroster
If we consider only the subspace $\Cal A_{\text{flat}}$ of 
flat $SU(3)$-connections, 
then the relevant strata are those of isotropy subgroup 
$\bold Z/3$, $U(1)$ and $SU(3)$, 
i.e. the irreducibles together with (2.1)(a) and (b).   
The reason is that, for our integral homology sphere $M$, 
there exist no nontrivial 
$U(1)$-representations $\pi_1(M)\to U(1)$ and hence every 
flat connection is gauge 
equivalent to the trivial connection.

Now, over $\Cal A$ there is the Chern-Simons functional 
$cs: \Cal A\to\bold R$ given by
$$
cs(A) = \frac{1}{8\pi^{2}} \int_X tr(A\wedge dA 
+ \frac{2}{3}A\wedge A\wedge A). \tag2.2
$$
With respect to a gauge transformation $g\in \Cal G$, we have
$$
cs(g\cdot A) = cs(A) + \deg (g) \tag2.3
$$
where $\deg g$ is the image, under $g^* = H^3(SU(3))\to H^3(X)=\bold Z$, of aanonical generator in $H^3(SU(3))$.  Because of (2.3), 
there is an induced mapping
$$
cs: \Cal B\to\bold R/\bold Z
$$
on the quotient spaces.  As is well-known \cite{T}, 
the gradient of $cs$ is given by
$$
\nabla cs(A) = -\frac{1}{4\pi^{2}} * FA, \tag2.4
$$
and so the set of critical points of $cs$ coincides with the moduli space
$$\align
\Cal M_{SU(3)} (X) &= \Cal A_{\text{flat}}/\Cal G \tag2.5 \\
&=\{[A]\in\Cal B\mid *F_A = 0\}
\endalign
$$
of gauge equivalent classes of flat $SU(3)$-connections on $X$.

By taking the intersection with the strata on 
$\Cal B$, we obtain an induced stratification on $\Cal M_{SU(3)}(X)$.  
In fact, because of (2.1), we can give an explicit 
description of all these strata.  First of all, we have the top strtum 
of irreducible, flat, $SU(3)$-connections denoted 
by $\Cal M^*_{SU(3)}$.  Then we have
the stratum consisting of $SU(3)$-connections which are the 
sum of an irreducible, flat, $SU(2)$-connection and a trivial product, 
$U(1)$-connection.  
Since this last stratum is isomorphic to the moduli 
space of irreducible, flat, $SU(2)$-connections, we wil 
denote it by $\Cal M^*_{SU(2)}$.  Finally, there is the 
stratum $[\theta]$ consisting of the single, isolated, trivial 
$SU(3)$-connection.

To obtain a well-defined invariant, Boden and Herald perturb the Chern-Simons 
functional so that the resulting critical points are finite number of 
regular points, i.e. points cut out transversely by the equation \cite{BH}.
Following the idea of Floer and others \cite{F} in $SU(2)$-gauge theory, they
consider the space $\Cal F$ of admissible perturbations consisting of
a collection of $n$ solid tori $\gamma_i:S^1\times D^2\to X, 1\leq i\leq n$,
and invariant functions $\tau_i:SU(3)\to \Bbb R$ and compactly supported 
2-form $\eta$ on $D^2$ with $\int_{D^2}\eta =1$. Then, for each element in 
$\Cal F$, the perturbation is given by adding to the Chern-Simons 
functional the following:
$$
h(A) = \sum_{i=1}^n\int_{D^2}\tau_i(hol_i(x,A))\eta(x)dx
$$
where $hol_i(x,A)$ is the holonomy of the connection $A$ around the loop 
$\gamma_i(S^1\times x)$.

Note that $h$ is invariant under gauge transformation and so $A\to cs(A)+h(A)$
descends to a function on $\Cal B$. After taking the differential, we obtain
a section of $\Cal A\times \Omega^1(X;su(3))$
$$
\align
\zeta_h: \Cal A &\longrightarrow \Omega^1(X;su(3)) \\
A &\longrightarrow \frac{-1}{4\pi^2}*F_A + \nabla h
\endalign
$$
A connection is said to be $h$-perturbed flat if it satisfies the equation
$\frac{-1}{4\pi^2}*F_A + \nabla h = 0$. The set of all gauge equivalent 
classes of such connections forms a moduli space, called the perturbed 
moduli space $\Cal M_{SU(3),h}(X) = \zeta_h^{-1}(0)/\Cal G$.
and has many properties
of $\Cal M_{SU(3)}$: For example it is compact (Proposition 2.9 of \cite{BH}).
In Theorem 3.13 of \cite{BH}, it is shown that, inside the space 
$\Cal F(\epsilon_0)$ of small ($\|h\|\leq\epsilon_0$), 
admissable perturbations, there exists
a Baire set $\Cal F(\epsilon_0)'$ of perturbations under which 
$\Cal M_{SU(3),h}(X)$ is regular. 
Moreover, for any two perturbations $h_{-1},h_1$,
in  $\Cal F(\epsilon_0)'$, there exits a path $h_t$ of 
small perturbations 
connecting $h_{-1},h_1$ such that the parametrized moduli space
$W=\{(A,t)\in \Cal A\times [-1,1]\mid \zeta_{h_t}(A)=0\}$ is also
regular.

The precise definition of small perturbation $h\in \Cal F(\epsilon_0)$ 
is in Proposition 3.7 of \cite{BH}. Basically, $\epsilon_0$
is chosen so that 

\roster
\item"(2.6)" If $\| h\|\leq\epsilon_0$ and $A$ is $h$-perturbed flat,
then there exists $\Hat A\in \Cal A_{flat}$ 
with $\| A-{\Hat A}_0\|\leq\epsilon_0$.
\item"(2.7)" If $A,A'$ are flat and lie in different components of
$\Cal A_{flat}$ then $\| A-A'\| \geq 2\epsilon_0$.
\endroster

\flushpar Since $\Cal A_{flat}$ is disjoint from those strata with 
isotropy subgroups $S(U(1)\times U(1)\times U(1)), S(U(1)\times U(2))$
we can choose $\epsilon_0$ so small that by (2.6) the 
perturbed moduli space 
$$
{\Cal M}_{SU(3),h}(X)={\Cal M}^{*}_{SU(3),h}(X)
\cup{\Cal M}^{*}_{S(U(1)\times U(1))}(X)\cup [\theta],
$$
in other words, a $h$-perturbed flat $SU(3)$-connectin is either irreducible
or with isotropy subgroup $U(1)$ or $SU(3)$.

Another consequence of (2.6), (2.7) is that associated to a $h$-perturbed
flat connection $A$, there is a unique component $\Hat A$ of flat connections
which is within $\epsilon_0$-distance. From this there is a well
defined invariant
$$
SF(\theta,A,h;s(u(1)\times u(2)))-2cs(\Hat A) 
$$
where the ambiguity of the path-dependent spectral flow 
$SF(\theta,A,h;s(u(1)\times U(2)))$ is cancelled by the
corresponding choice in $cs(\Hat A)$, as explained in \S 1.

We will need  several closely related spectral flows whose definitions 
can all be traced back to the linearized operator of $\zeta _h$ :
$$
*d_{A,h}=*d_A - 4\pi^2\cdot{\text {Hess}}~~
h(A):\Omega^1(X;su(3))\to\Omega^1(X;su(3))
$$
where Hess $h(A)$ is the Hessian of $h$. In terms of $*d_{A,h}$, there is the
self adjoint, Fredholm operator $K(A,h;su(3))$ given by
$$
\align
K(A,h,su(3)): &(\Omega^0\oplus\Omega^1)(X;su(3))\longrightarrow
(\Omega^0\oplus\Omega^1)(X;su(3)) \tag2.8\\
&(\xi,a)\longrightarrow (*d_A a, d_a\xi + *d_{A,h}(a)
\endalign
$$
Similarly, for a connection $A\in \Cal A$ with isotropy subgroup $U(1)$,
the structure group of $A$
can be reduced to $S(U(1)\times U(2))$. Hence we can form
the operator $K(A,h,s(U(1)\times u(2))$ by taking the tensor product 
of the self-adjoint operator in (2.8) with the adjoint representation 
$s(u(1)\times u(2))$. 

All the above are real self-adjoint operators, and so when we discuss its 
spectral flow we count the number of real eigenspaces crossing a 
$(-\epsilon/-\epsilon)$-reference
line. However, for a $S(U(1)\times U(1))$-connection $A$, we also have the
complex operator 
$K(A,h;\Bbb C^2)$ obtained by coupling the self-adjoint operator 
with the regular representation $\Bbb C^2$ of $S(U(1)\times U(1))$. Following
the convention in [BHKK], the spectral flows for these operators are referred
to the number of complex eigenspaces crossing the $(-\epsilon/-\epsilon)$-
reference line.

In the background of all these, there is also the deformation complex:
$$
\align
\Omega^0(X;su(3))& \overset{d_A}\to\longrightarrow
\Omega^1(X;su(3))\overset{*d_{A,h}}\to\longrightarrow \\
\Omega^1(X;su(3))& \overset{d_A^{*}}\to\longrightarrow
\Omega^0(X;su(3))
\endalign
$$
associated to a $h$-perturbed flat, $SU(3)$-connection $A$. In [BH],
it is shown that this is a Fredholm, elliptic complex with
$H^0(X;su(3))=\text{Ker} d_A$ 
and $H^1_{(A,h)}(X;su(3))=\text{Ker}(*d_{(A,h)})/Imd_A$.
In particular, for a $h$-perturbed flat $SU(3)$-connection $A$, we have
$$
\text{Ker} K(A,h;su(3)) = H^0(X;su(3))\oplus H^1_{(A,h)}(X;su(3)),
$$
and when $A$ is irreducible $H^0(X;su(3))=0$ and the vanishing of the kernel
of $K(A,h,su(3))$ is the same as the vanishing of $H^1_{(A,h)}(X;su(3))$.

Given a path $\{A_t\mid 0\leq t\leq 1\}$ of connections from the trivial 
$SU(3)$-connection, denoted by $\Theta$, to the connection $A=A_1$, we have
the family of self-adjoint, Fredholm operators $K(A_t,h;su(2))$ and hence
its spectral flow $SF(\Theta, A,h;su(3))$. Although the latter depends 
on the choice of paths, it only enters into our discussion through the 
expression $(-1)^{SF(\Theta,A,h;su(3))}$ for the sign. Since the ambiguity 
due to the choice of paths of $SF(\Theta,A,h;su(3))$ is 12 (see Prop 4.3
of \cite{BH}), this last sign is well-defined. 

Similarly for a path $\{A_t\mid 0\leq t\leq 1\}$ of 
$S(U(1)\times U(2))$-connections from the trivial representation, 
here denoted by $\theta$, we
have the spectral flows $SF(\theta,A,h;s(u(1)\times u(2)))$ and
$SF(\theta,A,h;\Bbb C^2)$ for the two families of self-adjoint operators
$K(A_t,h;s(u(1)\times u(2)))$ and $K(A_t,h;\Bbb C^2)$. The ambiguities due
to the choice of paths for $K(A_t,h;s(u(1)\times u(2)))$ are 8 and for
$K(A_t,h;\Bbb C^2)$ are 2. Once again we suppress this dependence because
they come into our application either as 
$(-1)^{SF(\theta,A,h;s(u(1)\times u(2)))}$
or as $SF(\theta,A,h;s(u(1)\times u(2)))-cs(\Hat A)$. 
Here, in the second case,  the ambiguities 
have been compensated by the Chern-Simons term.

With a choice of small perturbation $h$ which makes $\Cal M_{SU(3),h}(X)$
regular and with the convention of spectral flows as explained above,
Boden and Herald define their invariant $\lambda_{SU(3)}(X)$ by the formula
(1.2). The following is their main theorem (Theorem 1 of [BH]). 

\proclaim{Theorem 2.11 } 
Suppose $X$ is an integral homology $3$-sphere.
For generic small perturbation $h, \Cal M^{*}_{SU(3),h}(X)$ and
$\Cal M^{*}_{SU(2),h}(X)$ are smooth, compact, 0-dimensional manifolds. 
Choose a representative $A$ for each orbit $[A]\in \Cal M^{*}_{SU(3),h}(X)$ 
and in case $[A]\in \Cal M^{*}_{SU(2),h}(X)$ choose also a flat connection 
$\Hat A$ close to $A$. Define $\lambda_{SU(3)}(X)$ as in (1.2). Then
for $h$ sufficiently small, $\lambda_{SU(3)}(X)$ is independent of $h$
and the Riemannian metric and hence is a well-defined
topological invariant of $X$.
\endproclaim

\subhead 3.  Correction term via Floer chain complex \endsubhead

Recall that the reason for  introducing the Chern-Simons term
$cs(\Hat A)$ is to make the expression
$[SF(\theta, A, h; \CC^2) - 2cs(\Hat A) ]$
well defined, independent of the choice of path.
However there are other devices which can achieve
the same goal.

\proclaim{Lemma (3.1) }
If we use the same path $\{A_t \mid 0 \le t \le 1\},
A_0 = \theta, A_1 =A$ in computing the spectral
flows $SF(\theta,A,h;\CC^2),SF(\theta, A,h;s(u(1)\times u(2)))$,
then the difference  $\[SF\(\theta,A,h;\CC^2\)-
(1/4)\(SF\(\theta, A,h;s\(u(1) \times u(2)\)\)\)\]$
is well-defined, independent of the choice of paths $\lf\{A_t \mid 0 \le t
\le 1\rt\}$.
\endproclaim

\demo{Proof}
The ambiguities in $SF\(\theta,A,h;\CC^2\)$ and
$SF\(\theta, A,h;s(u(1)\times u(2)\)$ are the result of the nontrivial
nature of the fundamental group of the gauge space 
$\pi_1\( \CB(S(U(1)\times U(2))\)
= \pi_0 \(\text{Map}(X,U(2)\) = \ZZ$. A straightforward
computation shows that they are
8 for $SF\(\theta, A,h;s(u(1)\times u(2))\)$ and
2 for $SF\(\theta, A,h;\CC^2\)$. Hence, they
cancel out in taking the  difference
$SF(\theta, A,h;\CC^2))-(1/4)(SF(\theta, A,h;s(u(1)\times u(2)))$.
\enddemo

In view of (3.1),we can replace
$\lambda''_{SU(3)}(X)$ in (2.10) by the expression:
$$
\align
\La''_{SU(3)}(X) =& \tag3.2 \\
\sum_{ [A] \in \CM^*_{S(U(1) \times U(2)),h}(X)}
& (-1)^{SF(\theta,A,h;s(u(1) \times u(2)))}
\Big[SF\(\theta,A,h;\CC^2\)  \\
 & \qquad -(1/4)SF(\theta,A,h;s(u(1)\times u(2))) +(5/8)\Big]
\endalign
$$
This has the advantage that we can free ourselves
from the restriction of using only small perturbations.

On the other hand without the assumption of small
perturbation a new phenomenon has occurred.
Namely, during a parametrized family of perturbations $h_t$
a pair of $h_t$-perturbed connections $A_t(1), A_t(2)$
from different components of $\CM^*_{S(U(1)\times U(2))}$
can annihilate each other, as in the birth-death point
situation in parametrized Morse theory. In fact, as we will
see such an annihilation will cause a jump in the sum \thetag{3.2}
 and to compensate for this we have to introduce a
tertiary correction term from the Floer chain complex.

From now on, we consider the space of admissible perturbations $h
\in \CF$ without the assumption of being small, i.e.\
$(2.6),(2.7)$. Note that the choice of Wilson's
loops $\gamma_i : S^1 \times D^2 \rightarrow X$ and the
invariant functions $\tau_i : SU(3) \rightarrow \RR$ are
the same as in those in Floer's work. In particular, when we
restrict to the stratum $\CA_{S(U(1)\times U(2))}$,
we obtain the analogue of Floer's theory. Namely, we have
a chain complex $FC_*(X,h)$ over
$\ZZ/2$, which has the elements of $\CM^*_{S(U(1)\times U(2)),h}(X)$
as generators and is indexed by the Floer degree. This Floer degree
for a $h$-perturbed  flat connection $A$ is given by
$SF\(K(A_t,h,s(u(1)\times u(2))\)$ mod 8 where $A_t$ is any path  of
connections from the trivial connection $\theta$ to $A$.

Hence associated to $h$, we have the integer
$$
\align
\text{Floer }(X,h) = \sum_{p=0}^7 (-1)^p\, \dim_{\ZZ/2}
\lf\{ \text{ image of }
d: FC_{p+1}(X,h) \to FC_p(X,h)\rt\}
\endalign
$$
where the chain complex is a slight extension of Floer's treatment for $SU(2)$
to $S(U(1) \times U(2))$.
The associated Floer homology is the same since by concentrating on small
perturbations near $\CA_{SU(2)}$, we can deform
$FC_*(X,h)$ back to the $SU(2)$ situation. Note that the integer
$\text{Floer }(X,h)$ is sensitive to the perturbation $h$ and is precisely 
a device
which can account for the birth-death points between different perturbations.
With the Floer correction term as explained above, the perturbative
$SU(3)$-Casson invariant $\Lambda_{SU(3)}(X)$ of an integral homology
3-sphere $X$ is defined by the formula (1.3).

\flushpar{\bf Remark (3.3)}
As we will see in \S4,
the reason for $(5/8)$ in the formula of
$\Lambda''_{SU(3)}(X)$ is a
normalization factor to make sure that our invariant
has the property:
$\Lambda_{SU(3)}(-X) = \Lambda_{SU(3)}(X)$.
From Definition (1.3) it is clear that $8\cdot\Lambda_{SU(3)}$
is an integer; however $4\cdot\Lambda_{SU(3)}$ is
already an integer
because $\Sigma (-1)^{SF(\theta,A,h;s(u(1) \times u(2)))}$
is divisible by 2.

\proclaim{Theorem 3.4}
The number $\Lambda_{SU(3)}(X)$ is independent of the
Riemannian metric on $X$ and the admissible perturbation
$h \in \CF$ with the property
that the h-perturbed flat connections have
isotropy group $\Bbb/3$ or $U(1)$, and hence gives a well-defined, topological
invariant of the integral homology 3-sphere $X$.
\endproclaim

\demo{Proof} 
For the most part, we follow the argument
of Boden and Herold in \cite{BH} in establishing the well-definedness of
$\lambda_{SU(3)}(X)$. First of all, as in Theorem 3.13 of \cite{BH}, there
exists a  Baire set $\CF'$ of admissible perturbations
(not necessarily small) such that for $h \in \CF'$, an $h$-perturbed
flat connection $A$ has isotropy subgroup $\ZZ_3$ (irreducible case)
or $U(1)$ (reducible case). In the irreducible case, $\text{Ker}(K(A,h;su(3)))
=0$ and in the reducible case $\text{Ker}(K(A,h;\CC^2))=
\text{Ker}(K(A,h;s(u(1)\times u(2))) = 0$. These are referred to as the
regularity conditions because under these conditions the moduli spaces
$\CM^*_{SU(2),h}(X)$ and $\CM^*_{S(U(1)\times U(2)),h}(X)$ are smooth,
0-dimensional oriented compact manifolds. In particular, they consist
of finitely many points (up to gauge equivalence) and using the data
associated to them we can compute the sum
$\Lambda_{SU(3)}(X) = \Lambda'_{SU(3)}(X) + \Lambda''_{SU(3)}(X)
-(1/4)\text{ Floer }(X,h)$
as in (1.3).
\enddemo

Now for two such perturbations $h_0,h_1$, we can
connect them up by
a path of admissible perturbations $\rho=\{h(t)\mid 0 \le t \le 1\}$
such that the parametrized moduli space $W_{\rho}$ of $h(t)$-perturbed flat
connections is regular. More precisely, $W_{\rho} = W_{\rho}^* \cup
W_{\rho}^r$ with $W_{\rho}^*$ a  space of irreducible
$SU(3)$-connections and $W_{\rho}^r$ a space of
$S(U(1) \times U(2))$-connections.
Both $W_{\rho}^*$ and $W_{\rho}^r$
are properly embedded, smooth, oriented 1-manifold with boundary
where the boundary of $W_{\rho}^*$ is
the union $\CM^*_{SU(3),h_0} \cup \CM^*_{SU(3),h_1} \cup F$
with $F$ is a finite set of points in $W_{\rho}^r$.
and the boundary of $W_{\rho}^r$
is  $\CM^*_{S(U(1)\times U(2)),h_0} \cup \CM^*_{S(U(1)\times U(2)),h_1}$.

Note that $W_{\rho}^r$ may contain circle components. However, the regularity
condition for parametrized family implies that
they are finitely in number because
each gives rise to critical points with respect
to the projection in t-direction and there are finitely many
such critical points.
Thus by partitioning $[0,1]$ into small intervals
$[t(i),t(i+1)], 0 = t(0) < t(1)<\cdots < t(n)=1$ in a suitable fashion,
we can
break down these circles as a union of arcs whose intersection with
the closure $\bar{W^*_{\rho}}$ lie in the interior of these arcs. Since
$\Lambda_{SU(3),h_1} - \Lambda_{SU(3),h_0} = \sum \Sb{i=0}\endSb ^{n-1}
[\Lambda_{SU(3),h_{t(i+1)}} - \Lambda_{SU(3),h_{t(i)}}] $ is additive, we can
concentrate on the parametrized
families over these small intervals $[t(i),t(i+1)]$.
In short, we can assume that no circle components exist in $W^r_{\rho}$.

In view of the above discussion, let $S(0,1)$ denote the union of
curves in
$W^r_{\rho}$ that pass from $t=0$ to $t=1$, $S(0,0)$ denote
those that pass from
$t=0$ to $t=0$, and $S(1,1)$ to form $t=1$ to $t=1$. To simplify
our notation, we
list them as parametrized curves:
$$
\aligned
S(0,1) &= \lf\{ \gamma(j,u)\mid 0 \le u \le 1,\; j=1,\cdots,N\rt\}\\
S(0,0)& = \lf\{ \gamma'(j',u)\mid 0 \le u \le 1,\; j'=1,\cdots,N'\rt\}\\
S(1,1) &= \lf\{ \gamma''(j'', u)\mid 0 \le u \le 1,\; j''=1,\cdots,
N''\rt\}
\endaligned
$$

As we move along a curve $\{\gamma(j,u)\mid 0 \le u \le 1\}$ in $S(0,1)$
Taubes \cite{T} shows that the ``tangential'' signs
$(-1)^{SF(\theta,A,h;s(u(1)\times u(2))}$ at the two ends agree. Denote
this common value by $s_{n(j)}
= s_{n(\gamma(j,0))} = s_{n(\gamma(j,1))}$. On the other
hand, by \cite{BH} there are precisely $s_{n(j)}[SF(K(\gamma(j,u),h;
\CC^2)\mid 0 \le u \le 1) ]$ many $h$-perturbed flat, irreducible
$SU(3)$ connections sinking into or emitting
from this curve, each of which
is counted with sign $(-1)^{SF(K(A,h;su(3))}$.
Hence we have
$\text{Sum}(01) = \sum \Sb{j=1}\endSb ^N s_{n(j)}[SF(K(
\gamma(j,u),h;\CC^2) \mid 0 \le u \le 1)].$

Similarly, for a curve $\gamma'(j',u)\mid 0 \le u \le 1$ in $S(0,0)$
it follows from \cite{T} that the ``tangential signs''
$(-1)^{SF(\theta,A,h;s(u(1)\times u(2))}$
disagree. So we orient the curve in such a way that it traces form sign
$-1$ to sign $+1$. Then, in \cite{BH}, it is shown that there are
$-[SF(K(\gamma'(j',u),h;\CC^2)\mid 0 \le u \le 1) ]$ many $h$-perturbed flat,
irreducible, $SU(3)$-connections sinking into (or emitting from if negative)
points on this curve, counted with the signs, $(-1)^{SF(K(A,h;su(3))}$.
In toto, they give
$\text{Sum}(00) = \sum \Sb{j'=1}\endSb ^{N'}
-\[SF(K( \gamma'(j',u),h;\CC^2) \mid \le u \le 1)\]$.

The analysis for a curve $\gamma''(j'',u)\mid 0 \le u \le 1$ in
$S(1,1)$ is the same. From \cite{T}, the tangential signs
at the two ends disagree and we orient the curve so that it
travels from $-1$ to $+1$. From \cite{BH}, during its
history, there are precisely
 $+[SF(K(\gamma''(j'',u),h;\CC^2)\mid 0 \le u \le 1) ]$ many
 $h$-perturbed flat,
irreducible, $SU(3)$-connections sinking into (or emitting from)
points on this curve, counted with their signs, $(-1)^{SF(K(A,h;su(3))}$.
They give the sum:
$\text{Sum}(11) = \sum \Sb{j''=1}\endSb ^{N''}
-[SF(K( \gamma''(j'',u),h;\CC^2) \mid 0 \le u \le 1)].$

Note that an irreducible $SU(3)$-connection in $\CM^*_{SU(3),h_0}(X)$ at
$t=0$ can either travel all the way to $\CM^*_{SU(3),h_1}(X)$
at $t=1$ or be destroyed (likewise created) along the paths
in $S(0,1),S(0,0),S(1,1)$. In the first case, by \cite{T},
the contribution of the two end points cancel each other
in the difference $\Lambda'_{SU(3),h_0} - \Lambda'_{SU(3),h_1} $
while in the second case it enters as a term in
$-\text{Sum}(01), \text{Sum}(00), - \text{Sum}(11)$
(respectively for points created).
Thus we have the formula
$$
\align
\Lambda'_{SU(3),h_0} - \Lambda'_{SU(3),h_1}
= -\text{Sum}(01) + \text{Sum}(00) - \text{Sum}(11) \tag3.5
\endalign
$$

To prove (3.4), we add the term $\Lambda''_{SU(3),h_0} -
\Lambda''_{SU(3),h_1}$
to the two sides of (3.5) to get:
$$
\align
&\[\Lambda'_{SU(3),h_0} -  \Lambda''_{SU(3),h_0}\]-
\[\Lambda'_{SU(3),h_1} - \Lambda''_{SU(3),h_1}\]  \tag3.6 \\
&\hskip .66in =  -\text{Sum}(01) + \text{Sum}(00) - \text{Sum}(11)
+\[\Lambda''_{SU(3),h_0} - \Lambda''_{SU(3),h_1}\]
\endalign
$$
The idea is to rewrite the right hand side so that it
can be identified with the difference of Floer correction terms.
Note that, for a path
$\lf\{ \ga(u)\mid 0\le u\le 1\rt\}$ of
$S\(U(1)\times U(2)\)$-connections, the difference of the two
spectral flows
$$
\align
[ SF\( \theta, \ga(1), h;\CC^2 \)
&-\frac{1}{4} SF\( \theta,\ga(1), h; s \( U(1)\times U(2)\)\) ]\\
- [ SF\( \theta,\ga(0), h; \CC^2 \)
&-\frac{1}{4} SF\( \theta, \ga(0), h; s\(U(1)\times U(2)\)\)]
\endalign
$$
can be simplified into
$$
SF\[K(\ga(u),h; \CC^2)\mid 0\le u\le 1\]-
\frac{1}{4}SF\[ K(\ga(u),h; s(u(1) \times u(2)))
\mid 0\le u\le 1\]
$$
by the additivity of spectral flows.  We will apply
this device to the terms in $\La''_{SU(3), h_0} (X)
-\La''_{SU(3), h_1}(X)$
which correspond to pairs of points, connected up by paths
in $S(01), S(00), S(11)$.

For example, along a curve $\ga(j,u)$ in $S(01)$ the signs
$s_{(\ga(j,u))}$, at the two ends $u=0,1$ are the same, and so in the
difference $\La''_{SU(3), h_0}(X)-\La''_{SU(3), h_1}(X)$ we have
$$
\aligned
&s_{(\ga(j,1))}
\[ SF\( \theta, \ga(j,1),\, h;\CC^2\)
-\frac{1}{4} \,SF\( \theta, \ga(j,1),\,h; s(u(1)\times u(2))\) \]\\
-&s_{(\ga(j,0))} \[ SF \( \theta, \ga(j,0), h;\CC^2\)
-\frac{1}{4}
\,SF \(\theta, \ga(j,0),\,h; s (u(1) \times u(2) )\)\]\\
=& s_{(\ga(j,0))} \Big[ SF \( K \(\ga(j,u),\,h;\CC\) \mid 0\le u\le 1\)\\
& \qquad \qquad -\frac{1}{4}\, SF\( K \lf(\ga(j,u),\,h;s(u(1)\times u(2))
\mid 0\le u\le 1\rt)\)\Big].
\endaligned
$$
Note that the first sum cancels the contribution to the sum $S(01)$
by the same curve $\ga(j,u)$.

Similarly, along a curve $\ga'(j,u)$ in $S(00)$, we have the
following contribution to $\La''_{SU(3),h_0}(X)-\La''_{SU(3), h_1}(X)$:
$$
\aligned
&-s_{( \ga' (j', 1))} [SF \(\theta, \ga'(j', 1),h;\CC^2\)
-\frac{1}{4} SF( \theta, \ga'(j',1),h; s ( u(1)\times u(2))) ]\\
&-s_{( \ga'(j',0))} \[ SF \(\theta, \ga'(j',0),h;\CC^2\)
-\frac{1}{4} SF(\theta, \ga' (j',0),h;s( u(1)\times u(2))) \]\\
=&- \Big[ SF \( K \(\ga'(j',u),\,h;\CC^2\) \mid 0\le u \le 1\)\\
&\qquad - \frac{1}{4} \,SF\( K \(\ga'(j', u),\, h, s(u(1) \times u(2))\)\mid
0\le u\le 1\)\Big].
\endaligned
$$
In the last line, the first term cancels the corresponding contribution
to $\text{Sum}\,(00)$ in \thetag{3.6} by the curve.  The same works for a
curve $\ga''(j'',u)$ in $S(11)$ and provides us with the
contribution to $\La''_{SU(3),h_0}(X)-\La''_{SU(3),h_1}(X)$:
$$
\aligned
&s_{(\ga''(j'',1) )}
\[ SF\(\theta,\ga''(j'',1),h;\CC^2\)
-\frac{1}{4}SF\(\theta,\ga''(j'',1),h;s (u(1)\times u(2)) \) \]\\
+&s_{( \ga''(j'',0))}
\[SF( \theta, \ga''(j'',0),h;\CC^2)
-\frac{1}{4}SF(\theta,\ga'' (j'',0),h; s(u(1)\times u(2))) \] \\
=&+ \Big[ SF \( K  \(\ga''(j'', u),\, h;\CC^2 \) \mid 0\le u\le 1 \)\\
&\qquad \qquad - (1/4) \, SF
\( K  \( \ga''(j'', u),\,h ; s(u(1) \times u(2))\) \mid 0 \le u \le 1 \) \Big].
\endaligned
$$
Once again, this last term cancels the contribution
to $-\text{Sum}\,(11)$ in \thetag{3.6} by the same curve.

Thus we can rewrite \thetag{3.6} as follows:
$$
\align
&\[ \La'_{SU(3),h_0}(X) + \La''_{SU(3),h_0}(X)\]
- \[ \La'_{SU(3),h_0}(X)+ \La''_{SU(3),h_1}(X)\]
\\
=&\,(1/4) \[ - \text{Sum}' (01) + \text{Sum}'(00)-
\text{Sum}'(11)\].\tag3.7
\endalign
$$
Here the sums $\text{Sum}'(01), \text{Sum}'(00), \text{Sum}'(11)$
are obtained from the correspoinding sums $\text{Sum}(01),
\text{Sum}(00), \text{Sum}(11)$ by replacing the spectral flow
of the normal operator $K\( A_t, h; \CC^2\)$ by the corresponding
tangential operator
\newline
$K\(A_t, h; s(u(1)\times u(2)\)$ over the
same path of connections $A_t$.

To complete the proof of \thetag{3.4}, it remains to show that
the sum on the right hand side of \thetag{3.7} is $(1/4)\[
\text{Floer } (X, h_0)- \text{Floer }(X, h_1)\]$.  For this,
we observe that $\text{Sum}'(01) =0$ because by regularity the kernel
of the operator $K\( \ga(j,u),\,h; s(u(1)\times u(2)\)$ is zero
for every $u, \; 0\le u\le 1$.  On the
other hand, the spectral flows in $\text{Sum}'(00)$ and $\text{Sum}'(11)$
are not always zero as the kernels of $K\( \ga'(j',u),\,h; s(u(1) \times
u(2))\)$ and $K\( \ga''(j'', u),\,h, s(u(1)\times u(2)) \)$ may have
jumps at critical points of $t\(\ga'(j',u)\)$ and $t\(\ga''(j'',u)\)$.
The situation can be explained in terms of deformations of Floer
chain complexes.  In the language of parametrized Morse theory,
a Floer chain complex can be deformed from one to another by
a sequence of four moves:

\roster
\item"{\bf Move 1:}" (isotopy) The chain complex is unchanged

\item"{\bf Move 2:}" (handle slide) The chain groups are unchanged
but one of the differentials are changed by composing with an
elementary matrix
\item"{\bf Move 3:}" (birth point) Two new generators
$e,f$ are added in dimension $p,\, p+1$ with $de=f$, and
the rest of the chain complex is unchanged
\item"{\bf Move 4:}" (death point) the reverse.
\endroster
\smallskip

Furthermore, in the above Moves, the generators, other than those
pairs from birth-death points, move smoothly with constant Floer
index and zero tangential spectral flows.  While in a neighborhood
of a birth point in Move 3, we have pairs of generators with
consecutative Floer indice $p, p+1$. These pairs of generators trace
out a  curve $\{\ga(t), 0\le t \le 1\}$, and the tangential spectral flow
$SF (K (\ga(t), h; s(u(1)\times u(2))\mid 0\le t\le 1\}$ along this
curve equals 1 as it starts from index $p$ and ends at index $p+1$.
In the case of the death point, this is just the opposite.

Hence in Moves 1, 2, the expression $\[-\text{Sum}'(01)
+\text{Sum}'(00)-\text{Sum}'(11)\]$ is unchanged.  In Move 3,
this sum is increased by $(-1)^{p+1} \, (p+1)+(-1)^p\, p=(-1)^p$,
and in Move 4, it is decreased by $(-1)^p$. We now show that the
Floer correction term Floer $(X,h)$ changes in the same way.

Let $C_i(1), B_i(1), Z_i(1)$ be the $i^{\text{th}}$-chains,
$i^{\text{th}}$-boundaries, $i^{\text{th}}$-cycles associated
to the mod 2 Floer chain complex before making any move.
Let $C_i(2), B_i(2), Z_i(2)$ be the corresponding $\ZZ_2$-vector spaces
after one of the above moves. In Move 1, the dimension of all these
are unchanged since the Floer chain complexes before aand after are
identical.

For the second Move, the only changes are in the differentials
from $(p+1)$- to $p$-chains and from $p$- to $(p-1)$-chains, and so
$\dim B_i(1) =\dim B_i(2)$, for $i\neq p, \, p-1$.  As for $p,p-1$
terms, we have
$$
\aligned
\dim B_p(1) &=\dim C_{p+1}(1) - \dim Z_{p+1}(1)\\
&=\dim C_{p+1}(1) -\dim FH_{p+1}-\dim B_{p+1}(1).
\endaligned
$$
Since the last terms are the same for the chain complex after
the move, it follows that $\dim B_p(1)=\dim B_p(2)$.  Similarly,
we have $\dim B_{p-1}(1)=\dim Z_{p-1}(1) -\dim FH_{p-1}$.
As the latter are the same for both complexes, we have
$\dim B_{p-1} (1) =\dim B_{p-1}(2)$.  Consequently, in Move 2
the Floer correction term Floer $(X,h)$ is unchanged.

Consider the third Move where the dimension of $C_p(1),\,
C_{p+1}(1)$ are increased by $+1$ in going to $C_p(2),\,
C_{p+1}(2)$.  Again $\dim B_i(1)=\dim B_i(2)$ for $i\neq p+1,
p, p-1$.  As in the above but with degree shifting by 1, we have
$$
\aligned
\dim B_{p+1}(1) &= \dim C_{p+2}(1)-\dim Z_{p+2}(1)\\
&= \dim C_{p+2}(1) -\dim FH_{p+2} -\dim B_{p+2}(1) .
\endaligned
$$
Since these agree before and after, we have
$\dim B_{p+1}(1) =\dim B_{p+1}(2)$.  Using this last
equality, it also follows that
$$
\aligned
\dim B_p(1) &=\dim C_{p+1}(1) -\dim Z_{p+1}(1)\\
&=\dim C_{p+1}(1) - \dim FH_{p+1} -\dim B_{p+1}(1)\\
&=\[ \dim C_{p+1}(2)-1\]  - \dim FH_{p+1} - \dim B_{p+1}(2)\\
&=\dim B_{p+1}(2) - 1.
\endaligned
$$
Finally, by working from the lower degree end, we can deduce the formula
\newline $\dim B_{p-1}(1)
=\dim Z_{p-1} (1)-\dim FH_{p-1}$.  As these last terms are the same
for the chain complex after the Move, we have
$$
\dim B_{p-1}(1) = \dim B_{p-1}(2).
$$
Consequently, we can conclude that the Floer correction term
Floer $(X,h)$ is changed by $(-1)^p$ in Move 3.

Similarly, in Move 4, the Floer correction term Floer $(X,h)$ is
changed by $-(-1)^p$.  Since the argument is the same as above,
we will omit the details in here.   Thus we may conclude that for
a generic homotopy of perturbations the change in $\[-\text{Sum}'(01)
+ \text{Sum}'(00) -\text{Sum}'(11)\]$ is the same as the change in
$\text{Floer} (X, h)$.  This completes the proof that our invariant
$\La_{SU(3)}(X)$ is independent of all the choices.

\subhead 4. Properties of
$\BLa_{\bold S\bold U \boldkey( \bold 3 \boldkey)}
\boldkey(\bold X \boldkey)$ \endsubhead

The $SU(3)$-Casson invariants $\la_{SU(3)}(X)$ and $\La_{SU(3)}(X)$
are clearly different; nonetheless they share many common properties.
For example, if all the irreducible, flat,  $SU(2)$-connections
of $X$ are cut out tranversely, i.e.\ $H^1 \(X; s\(u(1)\times u(2)\)\)=
H^1 \(X, \CC^2\)=0$, then no perturbation along $\CM^*_{s(u(1)
\times u(2))}(X)$ stratum is necessary.  In this case, according
to Theorem 5.10 of [BHKK], the correction term $\la''_{SU(2)}(X)$
is given by
$$
\la''_{SU(3)}(X) =\sum_{[A] \in \CM^*_{SU(2)}(X)} (-1)^{SF
\(\theta, A, h; SU(2)\)} \[ \frac{1}{2}\,\rho \(K(A;\CC^2)\)\]
\tag 4.1
$$
where $\rho \(K (A;\CC^2)\)$ is the $\rho$-invariant of the self-dual
operator coupled to the regular representation of $SU(2)$.  A similar
result holds for $\La''_{SU(3)}(X)$.

\proclaim{Proposition (4.2)} Suppose $X$ is a homology $3$-sphere
with the property that every irreducible flat $SU(2)$-connection
$A$ has $H^1\(X; su(2)_A\)=0$ and $H^1\(X; \CC^2_A\)=0$. Then
there exist admissible perturbations $h$ which are zero on a
neighborhood of $\CM_{SU(2)}(X)$ and with respect to such perturbations:
$$
\aligned
&\La''_{SU(3)}(X) =\\
&\sum_{[A]\in \CM^*_{SU(2)}(X)} (-1)^{SF\(\theta,A;SU(2)\)}
\[\frac{1}{2} \,\rho\(K(A;\CC^2)\)-\frac{1}{8}\,\rho
\(K (A, su(2)\)\]
\endaligned
$$
\endproclaim
\demo{Proof}  
To calculate the spectral flows in $\La''_{SU(3)}(X)$,
we choose a path of connections $\lf\{A(t)\mid 0\le t\le 1\rt\}$ joining
the trivial connection $A(0)=\theta$ with an element $A(1)=A$
in the unperturbed moduli space $\CM^*_{SU(2)}(X)$.  Since $\CA_{SU(2)}$
is connected, we can choose the path lying inside $\CA_{SU(2)}$.
Note that along this path the coefficients $s\(u(1)\times u(2)\)$ is
decomposed into the sum $su(2)\oplus \RR$.  In particular, the kernel of
the
operator $K\(A(t);\RR\)$ from the second factor is constant and hence
gives no contribution to spectral flow, i.e. $SF\[ K\(A(t);\RR\)
\mid 0\le t\le 1\]=0$.
It follows that
$$
\aligned
&SF\[ K\(A(t); s\(u(1)\times u(2)\) \)\mid 0\le t\le 1\]\\
&= SF\[ K \(A(t); su(2)\) \mid 0\le t\le 1 \].
\endaligned
$$

From (5.4) and (6.5) of [BHKK], we have the following:
$$
\aligned
&SF[K( A(t);{\Bbb C}^2 ) \mid 0\leq t\leq 1] \\
&= 2 cs(A) +\frac{1}{2} \[\rho \(K \(A(1);\CC^2\) \)
-\rho \(K \(A(0); \CC^2\) \) \]\\
&\qquad +\frac{1}{2} \[ \dim \text{Ker} \(K \(A(1);\CC^2\)\)-\dim
\text{Ker} \(K\(A(0); \CC^2\)\)\],\\
&SF \[K\(A(t), su(2)\) \mid 0\le t\le 1\]\\
&= 8cs (A) +\frac{1}{2} \[\rho \(K \(A(1); su(2) \) \)
-\rho \(K \(A(0); su(2)\) \) \]\\
&\qquad +\frac{1}{2} \[\dim \text{Ker} (K( A(1); su(2) ))
-\dim \text{Ker}(K(A(0); su(2)))\].
\endaligned \tag4.3
$$
After substitution of \thetag{4.3} into $\La''_{SU(2)}(X)$, all
the terms except for the $\rho$-invariants cancel out
and the result is the formula in \thetag{4.2}.
\enddemo

\proclaim{Corollary (4.4)}
For the Brieskorn homology $3$-sphere $\Sigma(p,q,r)$, the difference
of the two $SU(3)$ Casson invariant
$(\la_{SU(3)} -\La_{SU(3)})(\Sigma(p,q,r))$ is given by
$$
\sum_{[A]\in \CM^*_{SU(2)}(X)} (-1)^{SF\(\theta,A;SU(2)\)}
\[\frac{1}{8}\rho\(K (A, su(2)\)\].
$$
\endproclaim

\demo{Proof}  
Note $\Sigma(p,q,r)$ satisfies the transversality
condition in \thetag{4.2}.  In addition, its Floer chain complex
is concentrated on odd degrees and so Floer $\(\Sigma(p,q,r),0\)=0$.
Our assertion follows immediately from comparing formulas in \thetag{4.4}
and \thetag{4.3}.
\enddemo

In (5.3) of [BH], it has been established that the invariant $\la_{SU(3)}(X)$
is independent of orientation.  We now show that
this is also true for the perturbative
$SU(3)$-Casson invariant  $\La_{SU(3)}(X)$.

\proclaim{Proposition (4.5)} $\La_{SU(3)}(-X)= \La_{SU(3)}(X)$.
\endproclaim
\demo{Proof}
We first consider the effect of reversing the orientation $X\to -X$
on the Floer chain complex $FC_*(X)$.  As in the usual Morse theory,
the effect of changing $X$ to $-X$ is accomplished by changing the
perturbed Chern-Simons functional by its negative and so
replaces $C_p(X)$ by its
dual $C^{-3-p}(X)=\text{Hom} \(C_{-3-p}(X),\ZZ/2\)$.  Thus,
we have
$$
\aligned
\dim B_p(-X) &=\dim \text{ Image } \[d: C_{p+1}
(-X)\to C_p(-X)\]\\
&= \dim \text{ Image } \[ d^*: C^{-3-(p+1)}(X)\to C^{-3-p}(X)\]\\
&=\dim \text{ Image } \[ d:C_{-3-p} (X) \to C_{-4-p}(X)\]\\
&=\dim B_{-4-p}(X),
\endaligned
$$
and so Floer $(X,h)=\text{Floer }(-X, -h)$.

On the other hand,
the spectral flows change via:
$$
\aligned
& SF_{-X} \(\Theta, A, -h; su(3)\) =-SF_X \(\Theta, A, h; su(3)\)-8\\
& SF_{-X} \(\theta, A, -h; s(u(1)\times u(2))\)
=-SF_X \(\theta, A, h; s(u(1)\times u(2)) \)-3\\
& SF_{-X} \(\theta, A, -h; \CC^2\) =-SF_X \(\theta, A, h; \CC^2\)-2
\endaligned
$$
Thus, changing the orientation leaves the signs of the
$SU(3)$-irreducibles $[A] \in \CM^*_{SU(3)}(X)$ unchanged as $(-1)^{-p-8}
=(-1)^p$.  On the other hand, for a $h$-perturbed flat,
$S\(U(1)\times U(2)\)$-connection $A\in \CM^*_{S\(U(1)\times U(2)\)}(X)$,
we have
$$
\aligned
&(-1)^{SF_{-X} \(\theta, A,-h;s\(u(1)\times u(2)\)\)}
\Big[SF_{-X} \(\theta, A,-h;\CC^2\) \\
&\qquad \qquad-\frac{1}{4} \(SF_{-X}
\(\theta,A,-h;\CC^2\)\)+\frac{5}{8}\Big] \\
=&-(-1)^{SF_X \(\theta,A,-h;s\(u(1)\times u(2)\)\)}
\Big[-SF_X \(\theta,A,-h;\CC^2\) -2  \\
&\qquad \qquad -\frac{1}{4}
\(-SF_X \(\theta,A,-h;s\(u(1)\times u(2)\) \)-3\)+\frac{5}{8} \Big] \\
=&-(-1)^{SF_X \( \theta,A,-h;s\(u(1)\times u(2)\) \)}
\Big[ SF_X \(\theta,A,-h;\CC^2\) \\
&\qquad \qquad -\frac{1}{4} \( SF_X \( \theta,A,-h;s\(u(1)\times u(2)\)\) \)
+\frac{5}{8}\Big].
\endaligned
$$
Consequently, our invariant $\La_{SU(3)}(X)$ is unchanged when we
reverse the orientation of $X$.
\enddemo

In [BH2], Boden and Herold showed that their $SU(3)$-Casson invariant
satisfy the connect sum formula:
$$
\align
\la_{SU(3)}\(X_1 \# X_2\) &= \la_{SU(3)} (X_1) + \la_{SU(3)}(X_2)\tag4.6\\
&+4\la_{SU(2)}(X_1) \cdot \la_{SU(2)}(X_2)
\endalign
$$
where $\la_{SU(2)}(X_i)$ is the normalized $SU(2)$-Casson invariant
(see [W]).  For the proof, they consider the connected sum $X_1\#
X_2$ as obtained from removing two flat $3$-balls $B_1, B_2$ from $X_1, X_2$
and gluing along the boundaries $X_1-B_1, X_2- B_2$ by an isometry.
Then they choose system of  loops in $X_1, X_2$ away from these balls
$B_1, B_2$, and based on these loops they choose admissable perturbations
$h_i$ of the self-dual equation on $\CA(X_i)$. The advantage for
this construction is that they can form the sum $h_1\# h_2$ perturbation
on $\CA(X)$ such that all the $h_1 \# h_2$-perturbation flat
connections are obtained from
gluing two $h_i$-perturbed
flat connections from $X_i$.  However, the  moduli
space $\CM_{SU(3), h_1 \# h_2}(X)$ obtained in this manner
is not necessarily regular.  Hence, they
have to choose an additional perturbation $h$ of $h_1 \# h_2$ to get
a regular moduli space $\CM_{SU(3), h}(X)$ for which they can compute
$\la_{SU(3)}(X)$ (see [BH2] for details).

To conclude this section, we obtain a similar connect sum
formula for $\La_{SU(3)}(X)$.

\proclaim{Theorem (4.7)} Let $X_1, X_2$ be integral homology $3$-spheres
and $X_1\# X_2$ be their connected sum.  Then,
$$
\aligned
&\La_{SU(3)} \(X_1 \# X_2\)
=\La_{SU(3)} \(X_1\)+ \La_{SU(3)} \(X_2\)+ \frac{9}{2} \,\La_{SU(2)}
\(X_1\) \La_{SU(2)}\(X_2\)  \\
&\qquad  - \frac{1}{4} \[ \text{Floer} \( X_1\# X_2, h\)
-\text{Floer} \(X_1, h_1 \) -\text{Floer} \(X_2, h_2\)\]
\endaligned
$$
where the perturbations $h_i$ for $X_i$ and $h$ for $X_1 \# X_3$
are the same as Boden-Herold perturbations in [BH2].
\endproclaim

\demo{Proof} 
As in \cite{BH2}, we choose small perturbations $h_1, h_2$
for the self-dual equations of $\Cal A_1, \Cal A_2$ such that
$\Cal M_{SU(3), h_i}^{*}(X_i) = \{A_{ij}\mid j=1,\cdots m_i\}$
and $\Cal M_{SU(2),h_i}^{*}(X_i) = \{B_{ij}\mid j=1,\cdots m_i\}$
consist of respectively isolated, $h_i$-perturbed flat
$SU(3)$-,$SU(2)$- connections. Then with respect to $h_1\#h_2$,
the perturbed flat connections in $\Cal A(X_1\#X_2)$ are given by
the glued connection $C_1\#C_2$ where $C_1,C_2$ ranges over the orbits
of $\{\theta_1, A_{1j}, B_{1k}\}\times \{\theta_2.A_{2j},B_{2k}\}$.
In particular, when the pair has isotropy subgroups $\Gamma_1,\Gamma_2$, 
then the glued
connections ranges over a connected
component isomorphic to
the double coset space $\Gamma_1\backslash SU(3)/\Gamma_2$.

As explained before, it requires a further perturbation $h$ to achieve
regularity. In \cite{BH2}, there is an explicit description of all the
resulting $h$-perturbed flat connections and their spectral flows as follows.

The pairs $A_{1j}\#\theta_2$ are single points and remain so after
$h$-perturbation. They are irreducible $SU(3)$-connections with
$$
SF_{X_1\#X_2}(\Theta,A_{1j}\#\theta_2,h; su(3))
= SF_{X_1}(\Theta,A_{1j},h;su(3)).
$$
The pairs $B_{1k}\#\theta_2$ are also single points and represent irreducible
$SU(2)$-connections with the same normal, tangential spectral flows as
the corresponding spectral flows of $B_{1k}$. In particular, the signed
coorection term for $B_{1k}$ in $\Lambda''_{SU(3)}(X_1\#X_2)$
is the same as the corresponding term for $B_{1k}$
in $\Lambda''_{SU(3)}(X_1)$. The same holds for the pair $\theta_1\#A_{2j},
\theta_1\#B_{2k}$. It follows that the contribution for these four type of
points to $\Lambda_{SU(3)}(X_1\#X_2)$ is the sum
$$
(\Lambda_{SU(3)}(X_1) + \frac{1}{4}{\text {Floer}}(X_1,h_1))+
(\Lambda_{SU(3)}(X_2) + \frac{1}{4}{\text {Floer}}(X_2,h_2))
$$

Next we consider the pairs $A_{1j}\#A_{2k}$, each of which yields
a component of $SU(3)$-irreducible connections isomorphic to $PSU(3)$.
Further perturbation by $h$ has the effect of introducing a Morse function $f$
to this component with its critical points $Q_{i,i'}$ as $h$-perturbed
flat connections associated to this component. The tangential spectral flow
$SF_{X_1\# X_2}(\Theta,Q_{i,i'},h;su(3))$ of $Q_{i.i'}$ is given by
$$
SF_{X_1}(\Theta_1,A_{1j},h_1;su(3)) + SF_{X_2}(\Theta_2,A_{2k},h_2;su(3))
+\text{ index of $f$ at}\, Q_{i,i'}
$$
Since we add up the signs $(-1)^{SF_{X_1\#X_2}(\Theta,Q_{i,i'},h,su(3))}$
in computing our invariant and since the Euler number of $PSU(3)$ is
zero, the total contribution of these points to our invariant
$\Lambda_{SU(3)}(X_1\#X_2)$ is zero.

In a similar manner, the pairs $A_{1j}\#B_{2k}$ yields a component
of $SU(3)$-connections isomorphic to $SU(3)/U(1)$. Since the Euler
number of the latter is zero, the same analysis shows that these pairs
give no contribution to $\Lambda_{SU(3)}(X_1\#X_2)$. Similarly for the pairs
$B_{1k}\#A_{2j}$, they again give no contribution.

There remain the pairs $B_{1k}\#B_{2k'}$, each of which gives rise to
a copy of
$U(1)\backslash SU(3)/U(1)$. However because the relative position of the two
$U(1)$'s there are two types of possible gluings with the result of
irreducible $SU(2)$'s and irreducible
$SU(3)$'s. In the situation of irreducible $SU(2)$'s, the double coset forms a
copy of $RP^3$ which, upon a Morse function perturbation, breaks into
four points $\{Q_{k,k',t},t=0,1,2,3\}$ indexed by the Morse index $t$.
The tangential spectral flow
$SF_{X_1\#X_2}(\theta,Q_{k,k',t},h;s(u(1)\times u(2)))$
of $Q_{k,k',t}$ is the sum $a_1+a_2+t$ where $a_1, a_2$ are respectively the
tangential spectral flows $SF_{X_1}(\theta_1,B_{1k},h_1;s(u(1)\times u(2))$,
$SF_{X_2}(\theta_2,B_{2k'},h_2;s(u(1)\times u(2))$ of $B_{1k},B_{1k'}$.
As for the normal spectral flows
$SF_{X_1\#X_2}(\theta, Q_{k,k',t},h;\Bbb C^2)$, they are the sum $b_1+b_2$
for all four points with $b_i$ the normal spectral flows
$SF_{X_i}(\theta_i,B_{ik},h_i;\Bbb C^2)$. Hence the normal contribution to
$\Lambda''_{SU(3)}(X_1\#X_2)$ by these four points is
$(-1)^{a_1+a_2}(1-1+1-1)=0$, or in other words the total contribution is zero.
As for the tangential contribution to $\Lambda''_{SU(3)}(X_1\#X_2)$, we have
$$
\align
&(-1/4)(-1)^{a_1+a_2}[(a_1+a_2)-(a_1+a_2+1)+(a_1+a_2+2)-(a_1+a_2+3)]\\
&=(1/2)(-1)^{a_1+a_2}.
\endalign
$$
Therefore, in toto the contribution of these irreducible $SU(2)$
representations is
$(1/2)\lambda_{SU(2)}(X_1)\cdot\lambda_{SU(2)}(X_2)$ as $\lambda_{SU(2)}(X_i)
=-\Sigma (-1)^{a_i}$.
Note that the constant term $(5/8)$ has no effect because it is counted with
the tangential signs and so gives $(1-1+1-1)=0$.

We still have to count the contribution from the $SU(3)$-irreducible
points in $B_{1k}\#B_{2k'}$. By making an equivariant Morse function
perturbation, each pair $B_{1k}\#B_{2k'}$ gives four irreducible $SU(3)$-
orbits $P_{k,k',t},t=0,1,2,3$ with identical sign $(-1)^{a_1+a_2}$. 
Consequently,
these four points give $4(-1)^{a_1+a_2}$ and the sum of all of them
is $4\lambda_{SU(2)}(X_1)\cdot\lambda_{SU(2)}(X_2)$.

Finally there are also changes in the Floer correction terms which are
compensated by
$$
\text{Floer}(X_1\#X_2,h) -\text{Floer}(X_1,h_1) - \text{Floer}(X_2,h_2).
$$
Adding up all these, we have the connect sum formula as
claimed.

\enddemo

\subhead \S5 Calculation of $SU(3)$-invariant for $(2,q)$-torus knots 
\endsubhead

Given a knot $T\subset S^3$, we have an integral homology $3$-sphere
$X(T,1/k)$ given by $1/k$-surgery of $T$. In turn, these homology 
$3$-spheres provide a sequence of $SU(3)$-invariants 
$\Lambda_{SU(3)}(X(T,1/k)), k=\pm 1,\pm 2 \cdots$ of the knot $T$. 
A natural question is the relation of these knot invariants to other
known knot invariants (c.f. (5.9)). 
As a first step, we consider in this section
the $(2,q)$-torus knot $T(2,q)$ and make explicit calcuclation of these
$SU(3)$-invariants.
\proclaim{Theorem 5.1}
Let $X_K$ denote the integral homology $3$-sphere given by $1/K$-
surgery of the $(2,q)$ torus knot $T(2,q)$. Then for $q=3,5,7,9,$
the $SU(3)$-Casson invariants $\Lambda_{SU(3)}(X_K)$ are as listed in the 
following table (5.2).
\endproclaim

\vskip2mm

\vbox{\tabskip=0pt \offinterlineskip
\def\tablerule{\noalign{\hrule}}
\halign to 400pt {\strut#& \vrule#\tabskip=1em plus2em&   
&#\hfil&\vrule#
&#\hfil&\vrule#
& #\hfil & \vrule#
\tabskip=0pt\cr\tablerule

&& (2,q)-torus knot
&& $\Lambda_{SU(3)}(X_K),K>0$
&& $\Lambda_{SU(3)}(X_K),K<0$ &\cr\tablerule

&& (2,3)
&& $\frac{1}{4}(10K-9)K$
&& $\frac{1}{4}(10K-11)K$ &\cr
\tablerule

&& (2,5)
&& $\frac{1}{4}(126K-79)K$
&& $\frac{1}{4}(126K-85)K$ & \cr
\tablerule

&& (2,7)
&& $\frac{1}{4}(540K-230)K$
&& $\frac{1}{4}(540K-242)K$ & \cr
\tablerule

&& (2,9)
&& $\frac{1}{4}(1540K-514)K$
&& $\frac{1}{4}(1540K-534)K$ & \cr
\tablerule}}

\centerline{\bf Table (5.2)}

\demo{Proof}
As is well-known, the $1/k$-surgery of a $(2,q)$
torus yields a Brieskorn sphere: $-\Sigma(2,q,2qk-1)$ for $K=k>0$ and
$\Sigma(2,2q,2qk+1)$ for $K=-k,k>0$ with the natural orientation from
singularity theory. The invariant $\lambda_{SU(3)}=\lambda'_{SU(3)}+
\lambda''_{SU(3)}$ of $\Sigma(2,q,2qk\pm 1)$ have been studied in great 
details in \cite{B},\cite{BHKK}. To simplify the notation, we write:

\flushpar {\bf (5.3)}

\roster
\item"(a)" $A(q,K)=\lambda'_{SU(3)}=\Sigma sign(B_j)$, sum over
the irreducible $SU(3)$-representations $B_j$ of $\pi_1(X_K)$.

\item"(b)" $B(q,K)=\lambda''_{SU(3)}=
(\epsilon/2)\Sigma \rho_{X_K}(A_j;\Bbb C^2)$, sum over the irreducible 
$SU(2)$-representations $A_j$ of $\pi_1(X_K)$.

\item"(c)" $C(q,K)=(-\epsilon/8)\Sigma \rho_{X_K}(A_j; su(2))$, 
sum over the irreducible
$SU(2)$-representations $A_j$ of $\pi_1(X_K)$.

\item"(d)" $D(q,K)={\text {Floer}}(X_K)$, the Floer correction term.
\endroster

Here $\epsilon$ is $-1$ for $K=k>0$ and $+1$ for $K=-k<0$. 
In terms of $A(q,K),B(q,K),C(q,K),D(q,K)$, we have
$$
\align
\lambda_{SU(3)}(X_K)&= A(q,K) + B(q,K)\\
\Lambda_{SU(3)}(X_K)&= A(q,K) + B(q,K) + C(q,K) + D(q,K).
\endalign
$$
In [B], Boden has shown that the irreducible $SU(3)$-representations $B_j$
of $\pi_1(X_K)$ all satisfy the regularity condition,
i.e. cut out transversely by equation,
and contribute with
$Sign(B_j)=(-1)^{SF(\Theta,B_j; su(3))}=1$. Thus no further perturbation is 
necessary, $h=0$, and $A(q,K)$ is the number of irreducible $SU(3)$-representations
of $\pi_1(X_K)$, listed in the first column of Table (5.4) below.

The aforementioned work of Boden can be regarded as an extension of the   
results on $SU(2)$-representations of $\pi_1(X_K)$, all of which satisfy 
the regularity condition. As in \cite{B2}, \cite{FS}, 
there are $(q^2-1)k/4$ of these $SU(2)$
representations $A_1,\cdots,A_{(q^2-1)k/4}$ which have odd spectral flow 
$SF(\theta,A_j;su(2))$ for $K=k>0$ and even for $K=-k<0$. It follows that 
the Floer chain complex has zero boundary map in all these cases. In 
particular, our Floer correction term $D(q,K)={\text {Floer}}(X_K)=0$
for all $K$.

In [BHKK], the terms $B(q,K)$ are computed and are listed in 
the second column
of (5.4). Since $A(q,K),B(q,K),D(q,K)$ are all known, our job is to 
calculate the remaining $C(q,K)$ in the third and fourth column in (5.4). Once
this is achieved, the proof of (5.1) is immediate by adding 
these columns together.

\vskip2mm

\centerline{\bf Table (5.4)}

\vskip4mm

\vbox{\tabskip=0pt \offinterlineskip
\def\tablerule{\noalign{\hrule}}
\halign to 423pt {\strut#& \vrule#\tabskip=1em plus2em&   
&#\hfil&\vrule#
&#\hfil&\vrule#
&#\hfil&\vrule#
& #\hfil & \vrule#
\tabskip=0pt\cr\tablerule

&& A(q,K)
&& B(q,K)
&& C(q,K) for $K > 0$
&& C(q,K) for $K < 0$ &\cr
\tablerule

&& && && && & \cr
&& $3K^2-K$
&& $\frac{K(-24K^2-84K+13)}{6(6K-1)}$
&& $\frac{K(12K^2+84K-11)}{12(6K-1)}$
&& $\frac{K(12K^248K-5)}{12(6K-1)}$ \cr
&& && && && & \cr
\tablerule

&& && && && & \cr
&& $33K^2-9K$
&& $\frac{K(-200K^2-1620K+151)}{10(10K-1)}$
&& $\frac{K(100K^2+1120K-87)}{20(10K-1)}$
&& $\frac{K(100K^2+48K-57)}{20(10K-1)}$ & \cr
&& && && && & \cr
\tablerule

&& && && && & \cr
&& $138K^2-26K$
&& $\frac{K(-784K^2-9128K+606)}{14(14K-1)}$
&& $\frac{K(392K^2+5992K-330)}{28(14K-1)}$
&& $\frac{K(392K^2+4816K-246)}{28(14K-1)}$ & \cr
&& && && &&  & \cr
\tablerule

&& && && && & \cr
&& $390K^2-58K$
&& $\frac{K(-2160K^2-33192K+1714)}{18(18K-1)}$
&& $\frac{K(1080K^2+20880K-890)}{36(18K-1)}$
&& $\frac{K(1080K^2+17640K-710)}{36(18K-1)}$ & \cr
&& && && && & \cr
\tablerule}}

\vskip4mm

Here the horizontal rows are the values of 
$A(q,K), B(q,K), C(q,K),K>0,C(q,K),K<0$ for q=3,5,7,9 respectively.
Recall that $C(q,K)$ is the sum of $\rho$-invariants of the adjoint 
representation $Ad(A_j)$ where $A_j$ runs through all the irreducible
$SU(2)$-representations of $\pi_1(\Sigma(2,q,2qk\pm 1))$. Our first step is
to tabulate these representations in a convenient manner. There are two cases,
$K>0$ and $K<0$, which have to be treated seperately.

\flushpar{\bf The Case $\Sigma (2,q,2qk-1), K>0$}

As is well-known, the Brieskorn sphere $\Sigma(2,q,2qk-1) = 
\Sigma(a_1,a_2,a_3)$, $a_1=2,a_2=q,a_3=2qk-1$
is a Seifert $3$-manifold with its Seifert invariant given by 
$(b_0,b_1,b_2,b_3)
=(-1,1,m,k), m=(q-1)/2$. As a Seifert manifold, 
its fundamental group has the following presentation:
$$
\align
&{\text{generators}}: x_1,x_2,x_3 \tag5.5 \\
&{\text {relations}}: x_1x_2x_3=h {\text{ central}}, x_1^2=h^{-1}, x_2^q=h^m, 
x_3^{2qk-1}=h^{-k}.
\endalign
$$

The central element $h$ plays an important role 
for an irreducible representation
$f:\pi_1(\Sigma(2,q,2qk-1))\to SU(2)$ because by Schur's lemma $f(h)=\pm 1$.
If $f(h)=1$, then $f(x_1)$ is also central and the representation becomes
abelian. As $H_1(\Sigma(2,q,2qk-1))=0$, this implies $f$ is the trivial 
representation. Hence we can omit this case and concentrate on $f(h)=-I$

Let $X_i=f(x_i)$. Then from (5.5) we have the following conditions:
$$  
X_1^2=-I, X_2^q=(-I)^m, X_3^{2qk-1}=(-I)^k, X_1X_2X_3=-I \tag5.6
$$
Consider an element $g\in SU(3)$ as a unit quarternion, written uniquely
in the form $g=cos\theta +sin\theta [icos(\pi t)+jsin(\pi t)], 
0\leq\theta <\pi$. Then the first
three conditions in (5.6) imply that 
$$
\align
trace(X_1)&=2cos(L_1/2\pi),\quad trace(X_2)=2cos(L_2/2\pi), \\
trace(X_3)&=2cos(L_3/2\pi)
\endalign
$$
where $L_1,L_2,L_3$ are integers with $L_1=1,0<L_2<q,0<L_3<(2qk-1), 
L_2 = m (mod2), L_3=k(mod(2)).$

In fact, by conjugation, we may assume that 
the pair $(X_1,X_2)$ takes the form $X_1=i$, 
$X_2 = cos(L_2\pi/q)+sin(L_2\pi/q)[icos(\pi t)+ jsin(\pi t)],$
with $0<t<1$. Such a choice of $(X_1,X_2)$ uniquely determines the 
representation because
$X_3=-X_2^{-1}X_1^{-1}
=icos(L_2\pi/q) + sin(L_2\pi/q) [cos(\pi t) + ksin(\pi t)]$.
Substition of this into $X_3^{2qk-1}=(-I)^k$
gives the constraint:
$sin(L_2\pi/2)cos(\pi t)$ $=cos(L_3\pi/2qk-1)$
on $t$. To solve this equation, we observe that as $t$ varies over $(0,1)$
the right hand side ranges monotonically over $(-1,1)$. Hence it is not 
difficult to work out the permissible values of $L_3$ for a fixed $L_2$.
For example, with $q=3$, we have $L_2=1$ and as $sin(\pi/3)=cos(\pi/6),
-sin(\pi/3)=cos(5\pi/6)$ the above constraint yields:
$\pi/6 < L_3\pi/(6k-1) < 5\pi/6$ and $L_3 = k$ mod $2$ which is
equivalent to $L_3 = (k-2)+2t, t=1,...,2k$. In this way we work out the 
following table of all admissible $L_1,L_2,L_3$ for $q=3,5,7,9$. 

\vskip4mm

\centerline{\bf Table 5.7}

\vskip4mm

\vbox{\tabskip=0pt \offinterlineskip
\def\tablerule{\noalign{\hrule}}
\halign to 400pt {\strut#& \vrule#\tabskip=1em plus2em&   
&#\hfil & \vrule #
&#\hfil&\vrule#
&#\hfil&\vrule#
&#\hfil&\vrule#
&#\hfil&\vrule#
& #\hfil & \vrule#
\tabskip=0pt\cr\tablerule

&& (2,q)
&& $L_1$
&& $L_2$
&& $L_3$
&&$t$
&&$e$ & \cr\tablerule

&& (2,3)
&& 1
&& 1
&& (k-2)+2t
&& 1...2k
&& 36k+12t-17 &\cr \tablerule

&& (2,5)
&& 1
&& 2
&& (k-2)+2t
&& 1...4k
&& 100k+20t-29 &\cr \tablerule

&& (2,5)
&& 1
&& 4
&& (3k-2)+2t
&& 1...2k
&& 160k+20t-33 &\cr \tablerule

&& (2,7)
&& 1
&& 1
&& (5k-2)+2t
&& 1...2k
&& 196k+28t-37 &\cr \tablerule

&& (2,7)
&& 1
&& 3
&& (k-2)+2t
&& 1...6k
&& 196k+28t-41 &\cr \tablerule

&& (2,7)
&& 1
&& 5
&& (3k-2)+2t
&& 1...4k
&& 280k+28t-45 &\cr \tablerule

&& (2,9)
&& 1
&& 2
&& (5k-2)+2t
&& 1...4k
&& 324k+36t-49 &\cr \tablerule

&& (2,9)
&& 1
&& 4
&& (k-2)+2t
&& 1...8k
&& 324k+36t-49 &\cr \tablerule

&& (2,9)
&& 1
&& 6
&& (3k-2)+2t
&& 1...6k
&& 432k+36t-57 &\cr \tablerule

&& (2,9)
&& 1
&& 8
&& (7k-2)+2t
&& 1...2k
&& 576k+36t-61 &\cr
\tablerule}}

\vskip2mm

\flushpar{\bf The case $\Sigma(2,q,2qk+1), K<0$} 

In this case, the Brieskorn sphere $\Sigma(2,q,2qk+1)=\Sigma(a_1,a_2,a_3),
a_1=2,a_2=q, a_3=2qk+1$ has its Seifert invariant given by 
$(b_0,b_1,b_2,b_3)=(1,-1,-m,-k)$. 
With these minor changes, the argument goes through the 
same way as before. We will omit the details and just
summarize our calculation in the following table.

\vskip2mm

\vbox{\tabskip=0pt \offinterlineskip
\def\tablerule{\noalign{\hrule}}
\halign to 400pt {\strut#& \vrule#\tabskip=1em plus2em&   
&#\hfil & \vrule #
&#\hfil&\vrule#
&#\hfil&\vrule#
&#\hfil&\vrule#
&#\hfil&\vrule#
& #\hfil & \vrule#
\tabskip=0pt\cr\tablerule

&& (2,q)
&& $L_1$
&& $L_2$
&& $L_3$
&&$t$
&&$e$ & \cr\tablerule

&& (2,3)
&& 1
&& 1
&& k+2t
&& 1...2k
&& 36k+12t+5 &\cr \tablerule

&& (2,5)
&& 1
&& 2
&& k+2t
&& 1...4k
&& 100k+20t+9 &\cr \tablerule

&& (2,5)
&& 1
&& 4
&& 3k+2t
&& 1...2k
&& 160k+20t+13 &\cr \tablerule

&& (2,7)
&& 1
&& 1
&& 5k+2t
&& 1...2k
&& 196k+28t+9 &\cr \tablerule

&& (2,7)
&& 1
&& 3
&& k+2t
&& 1...6k
&& 196k+28t+13 &\cr \tablerule

&& (2,7)
&& 1
&& 5
&& 3k+2t
&& 1...4k
&& 280k+28t+17 &\cr \tablerule

&& (2,9)
&& 1
&& 2
&& 5k+2t
&& 1...4k
&& 324k+36t+13 &\cr \tablerule

&& (2,9)
&& 1
&& 4
&& k+2t
&& 1...8k
&& 324k+36t+17 &\cr \tablerule

&& (2,9)
&& 1
&& 6
&& 3k+2t
&& 1...6k
&& 432k+36t+21 &\cr \tablerule

&& (2,9)
&& 1
&& 8
&& 7k+2t
&& 1...2k
&& 576k+36t+25 &\cr
\tablerule}}

\centerline{\bf Table(5.8)}

Now, as in \cite{FS}, the $\rho$-invariant $\rho_{X_K}(Ad(A_j))$ of the 
Adjoint representation can be computed by the following formula of
Dedekind sum:
$$
1/2\rho = 
3/2 + \sum_{i=1}^3\sum_{m=1}^{a_i-1}(2a/a_i)cot(\pi am/a_i^2)cot(\pi m/a_i)
sin^2(\pi em/a_i) \tag5.9
$$
where $e=\sum_{i=1}^3L_i(a/a_i)$ is listed in the last column of (5.7) (5.8).
With these data at hand, we can put them into (5.9) and then add up the 
$\rho$-invariants to get our formula for $C(q,K)$. In practice, this last step 
is a little easier. 
As in Lemma 10.3 of \cite{FS}, the sum 
$\Sigma (2/a_i)cot(\pi am/a_i^2)cot(\pi m/a_i)sin^2(\pi em/a_i)$
in (5.9)
is given by $int \Delta (x,y)-Area\Delta (x,y)$ where $\Delta (x,y)$
is the triangle with vertices $(0,0), (0,x), (x,y)$ and $Area\Delta (x,y)$
is its area and $int \Delta (x,y)$ is the number of the lattice points
inside and $(x,y)=(e,(b_i/a_i)e^{*}e)$. 
After putting in our data, we see that  
$\Delta (x,y)-Area\Delta (x,y)$ can be written 
a sum of greatest integer functions of the form 
[linear in i/linear in k] where i runs through integers in 
a fixed interval [0, const.k].
Then, after adding them up, $(2qk-1)C(q,K)$ 
can be shown to be a cubic polynomial in $K>0$
(respectively $K<0$). (This process is similar to the calculation of
$B(q,K)$ in \cite{BHKK}). Knowing that this is a cubic 
polynomial, the proof
reduces to
a simple matter of
linear algebra in deciding the coefficients by going through a finite
number of examples. In this way, we obtain 
the result as tabulated above and complete the proof of (5.1). 
\enddemo

\flushpar{\bf Remark(5.9)}
We conclude this paper with a conjecture. Observe that from our calculation
$\Lambda_{SU(3)}(X_K)$ are polynomials of degree 2:
$P_{+}(K,q)$ for $K>0$ and $P_{-}(K,q)$ for $K<0$. Moreover 
$P_{+}(K,q)=P_{-}(K,q) + (1/4)\vert K\vert N(q)$ 
where $N(3)=2, N(5)=6, N(7)=12,
N(9)= 20$. In all the cases computed here $N(q)\vert K\vert$ 
equals the number of 
irreducible $SU(2)$-representations. On the other hand, $\lambda_{SU(2)}$
is $K(q^2-1)/4$ where $-(q^2-1)/4$ is the second derivatives
$\Delta''_{T(2,q)}(1)$ at $+1$ of the normalized Alexander polynomial
of the $(2,q)$-torus knot $T(2,q)$. In view of this, a natural
{\bf conjecture} is that the SU(3)-knot invariants $\Lambda_{SU(3)}(X(T,1/K))$
are polynomials of degree 2 in $K$: $P_{+}(K)$ and $P_{-}(K)$ 
for $\vert K\vert$ large, and their difference are given by the
formula $P_{+}(K)=P_{-}(K) - \vert K \vert\cdot\Delta''_{T}(1)$.

\enddocument